%


%
%
\catcode`\@=11
%
\message{ *** SIAM Plain TeX macro package, version
   2.1.1, SEPTEMBER-1994 ***}
%
\font\tenrm=cmr10
\font\tenbf=cmbx10
\font\tenit=cmti10 
\font\tentt=cmtt10
\font\tensmc=cmcsc10
\def\tenpoint{%
   \def\rm{\fam0\tenrm}\def\bf{\fam\bffam\tenbf}%
   \def\it{\fam\itfam\tenit}\def\smc{\tensmc}\def\tt{\tentt}
        \textfont0=\tenrm \scriptfont0=\sevenrm
        \textfont1=\teni \scriptfont1=\seveni
        \textfont2=\tensy \scriptfont2=\sevensy
        \textfont3=\tenex \scriptfont3=\tenex
\baselineskip=12pt\rm}%

\font\ninerm=cmr9
\font\ninebf=cmbx9
\font\nineit=cmti9
\def\ninepoint{%
   \def\rm{\ninerm}\def\bf{\ninebf}%
   \def\it{\nineit}\baselineskip=11pt\rm}%

\font\eightrm=cmr8
\font\eightbf=cmbx8
\font\eightit=cmti8        
\font\eighti=cmmi8
\font\eightsy=cmsy8
\font\eightsmc=cmcsc10 at 8pt
\font\eighttt=cmtt8
   \def\eightpoint{%
   \def\rm{\fam0\eightrm}\def\bf{\fam\bffam\eightbf}%
  \def\it{\fam\itfam\eightit}\def\smc{\eightsmc}\def\tt{\eighttt}
        \textfont0=\eightrm \scriptfont0=\sixrm
        \textfont1=\eighti \scriptfont1=\sixi
        \textfont2=\eightsy \scriptfont2=\sixsy
\textfont3=\tenex \scriptfont3=\tenex
\baselineskip=10pt\rm}%

\font\sixrm=cmr6
\font\sixbf=cmbx6
\font\sixi=cmmi6        
\font\sixsmc=cmr5
\font\sixsy=cmsy6
\def\sixpoint{%
   \def\rm{\sixrm}\def\bf{\sixbf}%
   \def\smc{\sixsmc}\baselineskip=8pt\rm}%

\fontdimen13\tensy=2.6pt
\fontdimen14\tensy=2.6pt
\fontdimen15\tensy=2.6pt
\fontdimen16\tensy=1.2pt
\fontdimen17\tensy=1.2pt
\fontdimen18\tensy=1.2pt              

\def\rheadfont{\eightpoint\rm} 
\def\titlefont{\tenpoint\bf}
\def\authorfont{\eightpoint\smc}
\def\dedicatefont{\ninepoint\it}

\def\othersmcolon{;}
\catcode`\;= \active
\def;{\relax\ifmmode\othersmcolon 
           \else\ifdim\fontdimen1\the\font>0pt\/{\rm \othersmcolon}%
             \else\othersmcolon\fi\fi}

\hsize=31pc
\vsize=50pc
\parindent=2em

\newif\ifkeywords@
\newif\ifsubjclass@
\newif\ifdedicate@
\newif\ifrecdate@
\newif\ifoneclass@

\keywords@false
\subjclass@false
\dedicate@false
\recdate@false
\oneclass@false
     
\def\title#1\endtitle{\def\thetitle{\uppercase{#1}}%
          \def\\{\relax\ }\xdef\rightrh{\uppercase{#1}}}

\def\shorttitle#1{\xdef\rightrh{\uppercase{#1}}}

\let\protect\relax

\def\author#1\endauthor{%
        \def\and{\hbox{\sixrm AND }}\def\\{\break}
        \def\theauthor{\uppercase{#1}}%
         {%
           \def\\{\relax\ }
           \def\fnmark##1{}
           \def\and{and }
           \def\protect{\noexpand\protect\noexpand}
           \xdef\leftrh{\expandafter{\uppercase{#1}}}%
        }}

\gdef\fnmark#1{#1}
\gdef\address#1#2{\makefootnote@{\textfont2=\eightsy #1}{#2\unskip.}}
\def\journal#1{\def\thejournal{#1}}

\def\copyyear#1{\def\thecopyyear{#1}}
\copyyear{1989}
\def\vol#1{\def\thevol{#1}}
\def\no#1{\def\theno{#1}}
\def\date#1{\def\thedate{#1}}
\def\code#1{\def\thecode{#1}}
\def\dedicate#1{\dedicate@true\def\thededicate{#1}}
\def\keywords#1\endkeywords{\keywords@true\def\thekeywords{#1}}
\def\subjclass#1\endsubjclass{\subjclass@true\def\thesubjclass{#1}}
\def\oneclass{\oneclass@true}
\long\def\abstract#1{\def\theabstract{#1}}

\def\raggedcenter{\leftskip=0pt plus 1fill \rightskip=0pt plus 1fill}

\def\copyright#1{{\ooalign{\hfil\raise.07ex\hbox{c}\hfil\crcr#1\char"0D}}}

\def\recdate#1#2{\global\recdate@true
        \gdef\recdate@{#1Received by the editors \ignorespaces#2\unskip.}%
        \makefootnote@{}{\recdate@}}
                    
\let\topmatter=\relax
\def\endtopmatter{%
        \global\headline={\hss\vtop to \z@ {%
                \sixpoint\rm\noindent\thejournal \hfill%
                \rm\copyright{\sevensy}\rm%
                \ \thecopyyear\ Wiley\par
                \noindent Vol. \thevol , No. \theno , \thedate\hfill\thecode\par
        \vss}\hss}%
        \leavevmode\vskip8pt
        \vbox{\noindent\raggedcenter\let\\=\break\titlefont\thetitle
            \ifrecdate@*\fi}
        \vskip6pt
        \vbox{\noindent\raggedcenter\let\\=\break\authorfont\theauthor}
        \ifdedicate@
                \vskip6pt
                \vbox{\noindent\raggedcenter\let\\=\break\dedicatefont
                \thededicate}
        \fi
        \vskip14pt
        {\eightpoint{\bf Abstract.\ }\theabstract\par
        \ifkeywords@
                \vskip10pt
                {\bf Key words.} \thekeywords\par
        \fi
        \ifsubjclass@
                \vskip10pt
                {\bf AMS subject classification\ifoneclass@\else s\fi. }%
                \thesubjclass\par
        \fi}
        \vskip8pt
        \tenpoint}

\def\heading#1#2{%
        \vskip6pt{\bf #1.\enspace\ignorespaces#2.\enspace}\ignorespaces}

\def\thm#1{{\smc #1.}\begingroup\it\enspace\ignorespaces}
\let\lem=\thm
\let\cor=\thm
\let\prop=\thm

\def\endthm{\endgroup}
\let\endlem=\endthm
\let\endcor=\endthm
\let\endprop=\endthm

\def\prf#1{{\it #1}.\rm\enspace\ignorespaces}

\def\endproof{\vbox{\hrule height0.8pt\hbox{%
   \vrule height1.3ex width0.8pt\hskip0.8ex
   \vrule width0.8pt}\hrule height0.8pt
  }}



\def\ldisplaylinesno#1{\displ@y\halign{
  \hbox to\displaywidth{$\@lign\hfil\displaystyle##\hfil$}&
    \kern-\displaywidth\rlap{$##$}\kern\displaywidth\crcr
  #1\crcr}}

\def\meti#1{\parindent=2\parindent 
\par\indent\llap{#1\enspace}\ignorespaces\parindent=.5\parindent}



\newdimen\refindent@
\newdimen\refhangindent@
\newbox\refbox@
\setbox\refbox@=\hbox{\eightpoint\rm [00]}
\refindent@=\wd\refbox@

\def\resetrefindent#1{%
        \setbox\refbox@=\hbox{\eightpoint\rm [#1]}%
        \refindent@=\wd\refbox@}

\def\Refs{%
        \unskip\vskip2pc
        \centerline{\eightpoint\rm REFERENCES}%
        \penalty10000
        \vskip4pt
        \penalty10000
        \refhangindent@=\refindent@
        \global\advance\refhangindent@ by 2em
        \global\everypar{\hangindent\refhangindent@}%
        \parindent=0pt\eightpoint\rm}

\def\XRefs{%
        \unskip\vskip2pc
        \centerline{\eightpoint\rm REFERENCES}%
        \penalty10000
        \vskip4pt
        \penalty10000
        \refhangindent@=\refindent@
        \global\everypar{\hangindent\refhangindent@}%
        \parindent=0pt\eightpoint\rm}
                                      
\def\sameauthor{\leavevmode\vbox to 1ex{\vskip 0pt plus 100pt
    \hbox to 3em{\leaders\hrule\hfil}\vskip 0pt plus 300pt}}
    
\def\ref#1\\#2\endref{\leavevmode\hbox to \refindent@{\hfil[#1]}\enspace #2\par}

\def\xref\\#1\endref{\leavevmode #1\par}


\global\footline={\hss\eightpoint\rm\folio\hss}
\output{\plainoutput}
\def\plainoutput{\shipout\vbox{\makeheadline\pagebody\makefootline}%
  \advancepageno
  \ifnum\pageno>1
        \global\footline={\hfill}%
  \fi     
  \ifodd\pageno
        \global\headline={\hss\rightrh\hss{\tenpoint\rm\folio}}%
  \else
        \global\headline={\hskip-8pt{\tenpoint\rm\folio}\hss\leftrh\hss}%
  \fi
  \ifnum\outputpenalty>-\@MM \else\dosupereject\fi}
\def\pagebody{\vbox to\vsize{\boxmaxdepth\maxdepth \pagecontents}}
\def\makeheadline{\vbox to\z@{\vskip-22.5\p@
  \line{\vbox to8.5\p@{}\rheadfont\the\headline}\vss}%
    \nointerlineskip}
\def\makefootline{\baselineskip24\p@\vskip-8\p@\line{\the\footline}}
\def\dosupereject{\ifnum\insertpenalties>\z@ 
  \line{}\kern-\topskip\nobreak\vfill\supereject\fi}

\def\relaxnext@{\let\next\relax}
\def\footmarkform@#1{\ifmmode {}^{#1}\else$^{#1}$\fi }
\let\thefootnotemark\footmarkform@
\def\makefootnote@#1#2{\insert\footins
        {\interlinepenalty\interfootnotelinepenalty
        \eightpoint \splittopskip=\ht\strutbox
      \splitmaxdepth=\dp\strutbox
        \floatingpenalty=\@MM 
   \leftskip=\z@ \rightskip=\z@
        \spaceskip=\z@ \xspaceskip=\z@
        \leavevmode{#1}\footstrut\ignorespaces#2\unskip
        \lower\dp\strutbox\vbox to\dp\strutbox{}}}
\newcount\footmarkcount@
 \footmarkcount@=\z@                    
\def\footnotemark{\let\@sf=\empty \relaxnext@
        \ifhmode \edef\@sf{\spacefactor=\the\spacefactor}\/\fi
        \def\next@{\ifx[\next \let\next=\nextii@ \else
                \ifx"\next \let\next=\nextiii@ \else
                        \let\next=\nextiv@ \fi\fi\next}%
        \def\nextii@[##1]{\footmarkform@{##1}\@sf}%
        \def\nextiii@"##1"{{##1}\@sf}%
        \def\nextiv@{\global\advance\footmarkcount@\@ne
                \footmarkform@{\number\footmarkcount@}\@sf}%
        \futurelet\next\next@}
\def\footnotetext{\relaxnext@
        \def\next@{\ifx[\next \let\next=\nextii@ \else
                \ifx"\next \let\next=\nextiii@ \else
                        \let\next=\nextiv@ \fi\fi\next}%
        \def\nextii@[##1]##2{\makefootnote@{\footmarkform@{##1}}{##2}}%
        \def\nextiii@"##1"##2{\makefootnote@{##1}{##2}}%
        \def\nextiv@##1{\makefootnote@{\footmarkform@{%
                        \number\footmarkcount@}}{##1}}%
        \futurelet\next\next@}
\def\footnote{\let\@sf=\empty \relaxnext@
        \ifhmode \edef\@sf{\spacefactor\the\spacefactor}\/\fi
        \def\next@{\ifx[\next \let\next=\nextii@ \else
                \ifx"\next \let\next=\nextiii@ \else
                        \let\next=\nextiv@ \fi\fi\next}%
        \def\nextii@[##1]##2{\footnotemark[##1]\footnotetext[##1]{##2}}%
        \def\nextiii@"##1"##2{\footnotemark"##1"\footnotetext"##1"{##2}}%
        \def\nextiv@##1{\footnotemark\footnotetext{##1}}%
        \futurelet\next\next@}
\def\adjustfootnotemark#1{\advance\footmarkcount@#1\relax}

\skip\footins=18\p@ plus6\p@ minus6\p@

\def\footnoterule{\kern -4\p@\hrule width 3pc \kern 3.6\p@ } 

%
\catcode`\@=12


\def\defeq{\mathrel{\mathop=^{\rm def}}}
 
 
\topmatter
\journal{N{\smc ETWORKS}}
\vol{0}
\no{0, pp.~000--000}
\date{January 0000}
\copyyear{0000}
\code{000}

\title THE RING GROOMING PROBLEM\endtitle
 
\shorttitle{RING GROOMING}
 
 
\author Timothy Y. Chow\fnmark{$^{\dag}$} \and Philip J.
Lin\fnmark{$^{\ddag}$}\endauthor
 
\address{$^{\dag}$}{250 Whitwell Street \#2,
Quincy, MA 02169
({\tt tchow@alum.mit.edu})}
 
\address{$^{\ddag}$}{11 June Lane,
Newton, MA 02459
({\tt pjlin@ieee.org})}
 
\abstract{The problem of minimizing the number of bidirectional SONET rings
required to support a given traffic demand
has been studied by several researchers.
Here we study the related {\it ring grooming\/} problem
of minimizing the number of
add/drop locations instead of the number of rings;
in a number of situations this is a better approximation to
the true equipment cost.
Our main result is a new lower bound for the case of uniform traffic.
This allows us to prove that
a certain simple algorithm for uniform traffic
is in fact a constant-factor approximation algorithm,
and it also demonstrates that known lower bounds for
the general problem---in particular, the linear programming relaxation---are
not within a constant factor of the optimum.
We also show that our results for uniform traffic
extend readily to the more practically important case of
quasi-uniform traffic.
Finally, we show that if the number of nodes on the ring is fixed,
then ring grooming is solvable in polynomial time;
however, whether ring grooming is fixed-parameter tractable
is still an open question.}
 
\keywords ADM minimization, ring loading, BLSR, SONET over WDM,
 traffic grooming, parameterized complexity\endkeywords
 
\subjclass 90B10, 68W25, 94C30\endsubjclass
 
 
\endtopmatter
 
\heading{1}{Introduction}
Many of today's telecommunication networks
carry traffic on a configuration known as a {\it SONET BLSR\/}
or a {\it SONET bidirectional ring}.
A number of authors [2, 5, 11, 13, 17, 18]
have studied the so-called {\it ring loading problem},
which is the problem of minimizing the number of
concentric (or ``stacked'') SONET bidirectional rings required
to carry a given amount of traffic.
Some good approximation algorithms have been obtained for this problem.

The ring loading problem has the advantage of being easy to state
and amenable to rigorous analysis, but has the drawback that its cost function
gives only an approximation to the true cost of building a SONET ring network.
In some situations, particularly in so-called ``SONET over WDM'' networks,
a better approximation to the cost may be obtained by
defining the cost of a SONET ring to be proportional to
the number of its {\it add/drop points},
i.e., locations on the ring where traffic terminates
(i.e., begins or ends).
The idea is that at add/drop points, costly electronic multiplexers (ADM's)
must be installed, whereas elsewhere one can have
an ``optical passthrough'' or ``glassthrough''
whose cost is relatively negligible.

The problem of minimizing ADM's in SONET rings,
which is a special case of a more general problem known as
{\it traffic grooming},
has also attracted the attention of numerous researchers
[7, 12, 14, 15, 16, 19, 20, 21, 22, 23],
but these papers have concentrated either on
heuristics or ILP methods without any provable a priori performance bounds,
or on proving bounds under special assumptions
such as shortest-path routing.
(In the case of unidirectional rings, however,
there has been some recent work on approximation algorithms [1, 9].)
Part of the reason for this dearth of theoretical results is that
the ring grooming problem appears to be
more difficult to analyze than the ring loading problem.
Even in the case of uniform traffic on a bidirectional ring,
the answer to the following question was not known
prior to the present paper:
Does there exist an $\alpha$-approximation algorithm
for some absolute constant~$\alpha$?
Theorem~2 below answers this question affirmatively.
In fact, we show that if the traffic is {\it $K$-quasi-uniform\/}
(i.e., the largest traffic demand between any two nodes
is at most $K$ times the traffic between any other two nodes),
then there is an $\alpha$-approximation algorithm
with $\alpha$ depending only on~$K$.

We should add that our approximation algorithm
is not of tremendous practical value,
as its guaranteed approximation ratio is poor.
More important than the algorithm itself is Theorem~1,
a new lower bound that allows us to {\it prove\/} that our algorithm
is in fact a constant-factor approximation algorithm.
In certain cases, Theorem~1 is
an order of magnitude better than previous bounds
($O(n^{3/2})$ versus $O(n)$).
This has an important theoretical implication:
If one wishes to prove that a given algorithm---either
an existing algorithm or a new algorithm---is
a constant-factor approximation algorithm for the general case,
then one must first improve the best known lower bounds.

\heading{2}{Statement of problem}
We first give an informal description, which we follow with a precise
integer linear programming (ILP) formulation of the problem.

An instance of the ring grooming problem consists of
positive integers~$n\ge 2$ (the {\it ring size\/})
and~$c$ (the {\it capacity\/}),
and a list $\{j_1, k_1\}, \{j_2, k_2\}, \ldots, \{j_s, k_s\}$
of unordered pairs of integers between $1$ and~$n$ inclusive;
these are the {\it traffic demands.}
(For all~$i$, we require $j_i \ne k_i$,
but the list may contain the same pair of integers multiple times.)
Let $C_n$ be the cycle graph with $n$~vertices,
labeled 1 through~$n$ in clockwise order.
On this cycle are ``stacked'' multiple rings,
each with capacity~$c$ on each edge.
The traffic must be routed on these rings, i.e.,
for each traffic demand $\{j,k\}$ in the list,
we must choose one of the rings and
one of the two arcs on that ring between $j$ and~$k$
(either the ``minor arc'' or the ``major arc'') for its route.
To compute the cost of a routing,
we consider each of the stacked rings in turn.
If there exists any traffic terminating at
the $j$th vertex of the $i$th ring,
then an ADM must be placed on that ring at that vertex.
(Recall that ``terminating'' means either beginning or ending,
so in particular an ADM is required at both ends of a demand pair.)
The goal is to find a routing with the minimum number of~ADM's.

Now for the formal description.  Given an instance, we first set
$d_{jk}$ equal to the number of times the pair $\{j,k\}$
occurs in the list of traffic demands.
Note that $d_{jk} = d_{kj}$ and that $d_{jj} = 0$ for all~$j$.
We refer to $(d_{jk})$ as the {\it traffic matrix.}
We then set
$$ R = \sum_{j=1}^{n-1} \sum_{k=j+1}^n d_{jk}.$$
(This should be thought of as an upper bound on
the number of rings we will stack.)
To formalize major arcs and minor arcs,
we define the symbol~$\delta_{jkl}$ for $1\le j < k \le n$ and $1\le l\le n$
by setting $\delta_{jkl} = 1$ if
the arc in~$C_n$ between vertices~$j$ and~$k$
that contains the edge between vertices $l$ and~$l+1$
(if $l=n$ then this means the edge between vertices $n$ and~$1$)
also contains the edge between vertices $n$ and~$1$,
and setting $\delta_{jkl} = 0$ otherwise.

The variables of our ILP are of two types: 0-1 variables~$x_{ij}$
for $1\le i \le R$ and $1\le j\le n$
(indicating whether any traffic terminates at
the $j$th vertex of the $i$th ring)
and nonnegative integer variables $t^0_{ijk}$ and~$t^1_{ijk}$
for $1\le i \le R$ and $1\le j < k \le n$
(indicating the amount of traffic between vertices $j$ and~$k$ on ring~$i$
on each of the two arcs). 
The constraints are
$$\leqalignno{
\sum_{i=1}^R (t^0_{ijk} + t^1_{ijk})
  &= \hbox to 1in{$d_{jk}$,\hfil} \forall j, k \; (j<k)&(1)\cr
\sum_{j=1}^{n-1} \sum_{k=j+1}^n \biggl( \delta_{jkl} t^0_{ijk}
  + (1-\delta_{jkl}) t^1_{ijk} \biggr)
  &\le \hbox to 1in{$c$,\hfil} \forall i, l &(2)\cr
\sum_{k=1}^{j-1} (t^0_{ikj} + t^1_{ikj})
  + \sum_{k=j+1}^n (t^0_{ijk} + t^1_{ijk})
  &\le \hbox to 1in{$2cx_{ij}$,\hfil} \forall i, j&(3)\cr
0\le x_{ij} \le 1; \;\;\; t^0_{ijk} \ge 0; \;\;\;
  t^1_{ijk} &\ge \hbox to 1in{0\hfil} \forall i, j, k&(4)\cr
}$$
Constraint~(1) forces all the traffic to be routed,
constraint~(2) enforces the capacity constraint,
and constraint~(3) forces ADM's to be placed where traffic terminates.
The goal is to minimize the objective function $\sum_{i,j} x_{ij}$
subject to these constraints; formally,
$$\leqalignno{
m &:= \min\biggl\{ \sum_{i,j} x_{ij} \biggm| \hbox{$(x,t)$ satisfies
  constraints (1)--(4)}\biggr\}.&($\dag$)\cr
}$$

This completes the formal statement of the problem,
but some additional comments are in order.
Experts may notice that the version of the ring grooming problem
that we consider here differs in some details from
other versions in the literature.
These variations are discussed in section~7.

A trivial upper bound on~$m$ is~$2R$,
obtained by putting each unit of traffic on a separate ring.  From
this we see that if we choose any integer $R'>R$
and consider the ILP that is defined exactly as above
except with $R'$ in place of~$R$,
then the value of~$m$ will remain unchanged.
In other words, there is some leeway in the definition of~$R$;
any sufficiently large value will do,
but to be formal we need to fix a specific value.

The ring grooming problem is an optimization problem, and as usual 
it may be converted into a decision problem by introducing a bound~$M$
and asking if there exists a solution with $m\le M$.
Clearly this decision problem is in~NP.
Note, however, that if the traffic demands had been specified
not as a list but in the more succinct form of the traffic matrix
$(d_{jk})$, then it would no longer be clear that ring grooming
would be in~NP,
because then the number of rings required could potentially be
exponential in the size of the input.

In the rest of the paper we mostly use the informal language of
rings and ADM's rather than the formal ILP description.
So for example the reader should think of a variable number~$r$
of rings, all of which carry traffic,
rather than a fixed number~$R$ of rings, some of which may be~empty.

\heading{3}{Ring grooming compared to ring loading}
In view of the similarity between ring grooming and ring loading,
we might wonder whether the two problems are equivalent.
Note first that the ring loading problem as stated in~[18] does not
allow the demand between a given pair of nodes
to be split between the major and minor arcs,
whereas the ring grooming problem, at least as we have stated it, does.
Hence we should really compare ring grooming to
the special case of ring loading treated in~[5],
which allows traffic splitting and which can be solved in polynomial time.

As we shall see in the next section,
even with traffic splitting allowed,
ring grooming is NP-complete,
so the two problems are manifestly not equivalent.
To further highlight the difference,
consider the following instance.
Let $n=9$, let $c=1$, and let the list of traffic demands be
$$
\{1,2\}, \{1,3\}, \{2,3\}, 
\{4,5\}, \{4,6\}, \{5,6\}, 
\{7,8\}, \{7,9\}, \{8,9\}.
$$
One can check that the optimal ring grooming solution
uses three rings:
one with ADM's at vertices 1,~2, and~3,
one with ADM's at 4,~5, and~6,
and one with ADM's at 7,~8, and~9,
for a total of nine ADM's.
On the other hand, the optimal ring loading solution
uses two rings and fifteen ADM's:
one with ADM's at all nine vertices,
and one with ADM's at vertices 1, 3, 4, 6, 7, and~9.
In particular, the number of rings and the number of ADM's
cannot be simultaneously minimized.
This impossibility of simultaneous minimization
is what we would expect from
previous work on similar problems (e.g.,~[8]).

More complicated examples can easily be constructed,
and it appears that none of the techniques used successfully for ring loading
carry over readily to ring~grooming.

\heading{4}{Relationship to bin packing}
It is known~[15, 16] that traffic grooming on rings
is reducible from bin packing.
However, strictly speaking, the NP-completeness of
ring grooming as we have defined it
is not explicitly proved in the literature,
so we give a complete proof.

\prop{Proposition 1}
Ring grooming is NP-complete.
\endprop

\prf{Proof}
Recall that an instance of bin packing consists of
a positive integer~$B$ (the {\it bin size\/}) and
a set $A = \{a_1, \ldots, a_N\}$ of positive integers
(the sizes of the objects to be packed)
such that $a_i \le B$ for all~$i$.
The objective is to partition~$A$
into as few disjoint subsets as possible,
subject to the constraint that the sum of the~$a_i$
in each subset is at most~$B$.
Now, it is well known~[6] that bin packing is in fact
{\it strongly\/} NP-complete,
so we may (and will) assume that the $a_i$ are
given in {\it unary\/} rather than in {\it binary.}

Given an instance of bin packing, set $n=N+1$ and $c=B$.
For each $j$ in the range $1\le j \le n-1$,
add $2a_j$ copies of the pair $\{j,n\}$
to the list of traffic demands.
This gives us an instance of ring grooming
(whose size is polynomial in the size of the bin packing instance,
because the $a_j$ are given in unary).
If $m$ is the minimum number of ADM's for this instance,
then we claim that the minimum number of bins for the
original bin packing instance is $m-N$.
This implies in particular that
(the decision version of) ring grooming is NP-complete.

To prove the claim, we begin by observing that
we may assume that the optimal solution to the ring grooming instance
has the following property:
If $j < n$, then there is only one ring with an ADM at vertex~$j$.
For suppose we have a solution~$S$ with $r$ rings
in which there is some vertex $j<n$ at which
there is an ADM on rings $i_1, i_2, \ldots, i_s$ with $s\ge 2$.
Then we can create a new solution~$S'$ by adding a new ring to~$S$,
putting ADM's at vertices $j$ and~$n$ on this new ring,
transferring all the traffic between vertices $j$ and~$n$ from
rings $i_1, i_2, \ldots, i_s$ to the new ring,
and deleting the ADM's at vertex~$j$ from rings $i_1, i_2, \ldots, i_s$.
This is always feasible, because the total traffic
between vertices $j$ and~$n$ is~$2a_j$,
which is at most $2B = 2c$ and therefore fits onto a single ring
(we route $a_j$ units on the major arc and $a_j$ units on the minor arc).
Moreover, since $s\ge 2$, when we pass from $S$ to~$S'$
we add two ADM's and delete at least two ADM's,
so $S'$ contains no more ADM's than~$S$.

For the remainder of the proof we restrict ourselves to configurations
such that every ring has an ADM at vertex~$n$
(this is necessary because all traffic terminates at vertex~$n$)
and such that, for each $j<n$, exactly one ring has an ADM at vertex~$j$.
So if $r$ is the number of rings, then the number of ADM's is $m=r+N$.
Since $N$ is given, minimizing~$m$ is equivalent to minimizing~$r$.

All that remains is to establish a one-to-one correspondence
between rings and bins.
Specifically, we need to show
that for any $J \subseteq \{1, 2, \ldots, N\}$, we have
$$\sum_{j\in J} a_j \le B$$
if and only if all the demands $\{d_{jn} \mid j\in J\}$
can be routed on a ring with ADM's at vertex~$n$ and all vertices $j\in J$.
(As always, $(d_{jk})$ is the traffic matrix.)
Suppose first that $\sum_{j\in J} a_j \le B = c$.
Then for each $j\in J$,
we need to route $d_{jn} = 2a_j$ units of traffic between vertices $j$ and~$n$.
We do this by routing $a_j$~units on the major arc
and $a_j$~units on the minor arc;
the given inequality ensures that no edge's capacity is exceeded,
and therefore the proposed routing is feasible.
Conversely, given any feasible routing,
consider the ``cut'' that separates vertex~$n$ from the rest of the ring.
There are two edges across this cut, each with capacity~$c$,
so the total amount of traffic across the cut is at most~$2c$, i.e.,
$$2c \ge \sum_{j\in J} d_{jn} = \sum_{j\in J} 2a_j,$$
or $\sum_{j\in J} a_j \le c = B$ as required.\qquad\endproof

If the bin size is fixed, then bin packing
is solvable in polynomial time~[6].
A similar result holds for ring grooming.

\prop{Proposition 2}
If $n$ is fixed, then ring grooming is solvable in polynomial time.
\endprop

\prf{Proof}
Create a graph whose vertices are vectors with $n(n-1)/2$ coordinates,
with one coordinate corresponding to
each pair $\{j,k\}$ of distinct integers from $1$ to~$n$,
and with each coordinate taking on a value between zero and~$d_{jk}$.
Note that the number of such vectors is polynomial in the
size of the instance.
If $u$ and~$v$ are two such vectors,
then we draw an edge from $u$ to~$v$
if the coordinates of $v-u$ are all nonnegative
and if the traffic demands represented by $v-u$
all fit onto a single ring of capacity~$c$.
(Determining whether such an edge exists
requires solving a ring loading problem with traffic splitting,
which as we remarked before can be solved in polynomial time~[5].)
The weight of such an edge, if it exists,
equals the number of ADM's required to support the single ring
with the demands $v-u$.
Ring grooming now reduces to finding the minimum-weight path
from the zero vector to the vector with coordinates~$d_{jk}$.\qquad\endproof

Although theoretically polynomial-time,
the algorithm of Proposition~2 is not practical
for typical values of $n$ and~$c$ (i.e., $n\le 16$, and $c$ in the hundreds).
Intuitively, this is because the degree of the polynomial depends on~$n$;
in the terminology of parameterized complexity~[4],
Proposition~2 shows only that ring grooming is in~XP
and not necessarily in~FPT (if $n$ is the parameter).
Michael Fellows (personal communication) has shown that
if the number of ADM's is taken to be the parameter,
then ring grooming is in FPT.
Unfortunately, the number of ADM's tends to be
much larger than the ring size,
so it remains an interesting open problem
whether ring grooming is in FPT if $n$ is the parameter.

\heading{5}{Lower bounds}
The simplest lower bound on~$m$ is the {\it LP bound,}
obtained by relaxing the integer variables to rational (or real) variables.
It is easy to show that the LP bound is given by the explicit formula
$$ m \ge \sum_{j=1}^n \sum_{k=1}^n {d_{jk} \over 2c}.$$
Unfortunately, although this bound is easy to compute,
we shall see shortly that the integrality gap is unbounded
(this is one reason why ring grooming appears to be hard to solve).
A slightly better bound may be obtained by observing that $\sum_k d_{jk}$
is the total amount of traffic terminating at vertex~$j$
and that each ADM can handle at most $2c$ units of traffic.
Since the number of ADM's is an integer, this yields
$$m \ge \sum_{j=1}^n
  \biggl\lceil \sum_{k=1}^n {d_{jk} \over 2c} \biggr\rceil.$$
We call this the {\it add/drop lower bound}.
Notice that it is slightly better than
rounding up the LP bound to the nearest integer,
since we round up at each vertex~$j$ separately before summing over~$j$.

It is tempting to wonder if there is any sense in which
the ceiling signs may be pushed inside the inner sum as well.
This cannot be done na\"{\i}vely,
but in [7] the authors give a more delicate argument along these lines
that yields a lower bound that is
sometimes better than the add/drop lower bound.
To state this bound we need some notation.
Define $q_{jk}$ and~$r_{jk}$ to be
the quotient and remainder when $d_{jk}$ is divided by~$c$, i.e.,
choose $q_{jk}$ and~$r_{jk}$ such that
$$d_{jk} = cq_{jk} + r_{jk} \qquad \hbox{with $0 \le r_{jk} < c$.}$$
Next, order the $d_{jk}$ with $j<k$
in such a way that
their corresponding remainders $r_{jk}$ decrease monotonically.
For simplicity, we write
$D_p$, $Q_p$, and~$R_p$ respectively for $d_{jk}$, $q_{jk}$, and~$r_{jk}$,
where $p$ runs from $1$ to $n(n-1)/2$ and the labeling is chosen so that
$$\hbox{$R_p \ge R_{p'}$ whenever $p < p'$.}$$

\prop{Proposition 3 [7]}
Suppose we are given an instance of the ring grooming problem.
With the notation above, let $P$ be the smallest integer such that
$$\biggl( \sum_{p=1}^{n(n-1)/2} cQ_p \biggr)
  + \biggl( \sum_{p=1}^P {R_p + c \over 2} \biggr)
  \ge \sum_{p=1}^{n(n-1)/2} D_p.$$
Then regardless of the routing,
the minimum number~$m$ of ADM's as defined in $(\dag)$ must satisfy
$$m \ge P + \sum_{p=1}^{N(N-1)/2} Q_p.$$
\endprop

The add/drop lower bound and
the lower bound of Proposition~3
are the best known lower bounds in general.
However, for $K$-quasi-uniform traffic, a better lower bound is possible,
and this is the first main theorem of this paper.
Recall from the introduction that ``$K$-quasi-uniform'' means that
$d_{jk}/d_{j'k'} \le K$ for all $j\ne k, j'\ne k'$.
In particular, $d_{jk} \ne 0$ for all $j\ne k$.
Conversely, if $d_{jk} \ne 0$ for all $j\ne k$,
then the traffic is $K$-quasi-uniform for some $K\ge 1$.
Uniform traffic corresponds to the case $K=1$.
Since bidirectional rings are often used in the core of a network,
where every node typically has traffic to every other node,
quasi-uniform traffic is not a bad approximation to reality.
However, we shall see below that our results give bounds that depend on~$K$,
so we are still far from satisfactorily solving arbitrary instances
of ring grooming.

We need some further terminology and notation.
The {\it bandwidth\/} consumed by a traffic demand
is the capacity of the demand multiplied by
the number of edges that it traverses from end to end.
For instance, in the example given in section~3,
the solution using shortest-path routing
consumes 12~units of bandwidth,
while the nine-ADM solution consumes 27~units of bandwidth.

Our results for quasi-uniform traffic are
really just easy corollaries of the case of uniform traffic,
so let us now focus attention on uniform traffic.
Then all the $d_{jk}$ with $j\ne k$ are equal,
and we use the letter~$d$ to denote this common value.
It is also convenient to define
$$f \defeq {d \over 2c}.$$

\thm{Theorem 1}
With the notation above, the minimum number~$m$ of ADM's
required for uniform traffic must, regardless of routing, satisfy
$$m \ge {(n^2-1)\sqrt{f} \over 4}.$$
\endthm

\prf{Proof}
Let $S$ be any feasible solution to the given instance.
For each ADM~$A$ in the solution, let $R(A)$ denote the ring that $A$ is on,
and define $B(A)$ by
$$B(A) \defeq { \hbox{Total bandwidth of the traffic carried on $R(A)$} \over
   \hbox{Total number of ADM's on $R(A)$} }.$$
Now suppose that on a particular ring in~$S$, there are $x$ ADM's.
Then there are $x(x-1)/2$ pairs of ADM's,
and between each pair there are at most~$d$ units of traffic on this ring
(since there are $d$ units total of traffic between this pair).
The circumference of the ring is $n$~edges,
so the total bandwidth on this ring is at most $dnx(x-1)/2$.

On the other hand, the total bandwidth on this ring is at most $cn$
just by capacity constraints.  Therefore, if $R(A)$ has $x$~ADM's, then
$$\eqalignno{B(A) & \le {\min(dnx(x-1)/2, cn) \over x} = cn \cdot \min(f(x-1), 1/x)\cr
   &\le cn\cdot \min(fx, 1/x) \le cn \sqrt{f}.\cr}$$
(The last inequality in this chain may be obtained, for example,
by graphing the curves $y=fx$ and $y=1/x$ and
finding their point of intersection,
or by noting that the minimum of two positive reals
is less than or equal to their geometric mean.)
This inequality is independent of~$x$,
so if we sum over all ADM's~$A$ in~$S$
we find that the total bandwidth~$B$ used by all rings satisfies
$$B \le mcn\sqrt{f}.$$

To obtain a lower bound for~$B$,
note that shortest-path routing uses the minimum amount of bandwidth.
It is well known and easy to prove that
$$B \ge {dn(n^2-1) \over 8} = {fcn(n^2-1)\over 4}.$$
(If $n$ is even, one can improve this bound to $dn^3\!/8$,
but we ignore this.)
Combining the lower and upper bounds for~$B$
and solving for~$m$ yields the theorem.\qquad\endproof

To compare this bound with the add/drop lower bound
and the lower bound of Proposition~3,
consider the special case in which $n=2c+1$
(where $c$ as always is the capacity) and $d=1$.
Then the add/drop lower bound is~$n$
and the bound of Proposition~3 is about~$2n$,
while the lower bound of Theorem~1 is about ${1\over 4}n^{3/2}$.
Thus, for the small values of~$n$ that typically occur in practice,
the bound of Theorem~1 is actually worse than the previous bounds;
however, as explained in the introduction,
it is a significant improvement theoretically.

Note that the proof of Theorem~1 never uses
the integrality of the traffic variables $t^0_{ijk}$ and~$t^1_{ijk}$.
The lower bound therefore also holds for the ``semi-relaxation''
of ring grooming in which these variables (but not the~$x_{ij}$ variables)
are allowed to be arbitrary real numbers.

\cor{Corollary 1}
If an instance of ring grooming is $K$-quasi-uniform,
with $d$ being the largest value of the $d_{jk}$,
then the minimum number~$m$ of ADM's must satisfy
$$m \ge {(n^2-1)\sqrt{d/(2cK)} \over 4}.$$
\endcor

\prf{Proof}
Between any two nodes there are at least $d/K$ units of traffic,
and by Theorem~1, supporting just this subset of the traffic
already requires at least
${1\over 4} (n^2-1)\sqrt{d/(2cK)}$~ADM's.\qquad\endproof

\heading{6}{Approximation algorithm for $K$-quasi-uniform traffic}
Again, we focus first on the case of uniform traffic.
Given $n$, $c$, and~$d$, we define $f=d/2c$ as before.
Then our approximation algorithm, which we call Algorithm~A,
proceeds as follows.
If $f\ge 1$ then simply create, for each pair of vertices $j$ and~$k$,
a set of $\lceil f \rceil$ rings with just two ADM's---one at~$j$
and one at~$k$---and route all the traffic between $j$ and~$k$
on these rings (and then terminate the algorithm).
Otherwise, let
$$M = \cases{\bigl\lfloor \sqrt{2/f} \bigr\rfloor
      ,&if $2\le \bigl\lfloor \sqrt{2/f} \bigr\rfloor \le n$;\cr
  2,&if $\bigl\lfloor \sqrt{2/f}\bigr\rfloor < 2$;\cr
  n,&if $\bigl\lfloor \sqrt{2/f}\bigr\rfloor > n$.\cr}$$
Construct an {\it $(n,M,2)$-covering design},
i.e., a family of $M$-element subsets,
called {\it blocks}, of $\{1,2,\ldots,n\}$
such that every pair of integers $(j,k)$ with $1\le j < k \le n$
is contained in at least one of the blocks.
Finally, take each block $\{v_1, v_2, \ldots, v_M\}$ in turn,
and construct a ring with an ADM at each~$v_i$,
routing {\it all\/} the demands among these~$v_i$ on this ring
(provided that these demands have not already been routed on a previous ring).
If at the end of this process there are some ADM's
at which no traffic terminates, then discard them, and terminate.

This completes the description of Algorithm~A, 
but several steps require further elaboration and justification.
First we need a fast method of constructing an $(n,M,2)$-covering design.
This is easy: Let $\mu=\lfloor M/2\rfloor$,
and let $F$ be a family of $\mu$-element subsets of $\{1, 2, \ldots, n\}$
whose union is the entire set.
Assume further that $F$ has the minimum possible size,
i.e., that it has $\lceil n/\mu\rceil$ members.
We obtain an $(n,M,2)$-covering design of size at most
$$\lceil n/\mu\rceil \choose 2$$
by taking all pairs of members of~$F$
(if $M$ is odd, then we must
add an additional integer---it does not matter which one---to each such pair
in order to make the size of each block exactly equal to~$M$)
and discarding any duplicates.

Next, we need to show that $M$ is small enough
for all the designated traffic to fit onto a single ring,
and we also need to provide a fast method
for computing an actual routing.
This is accomplished by the following lemma and its proof.

\lem{Lemma 1}
Let $f=d/2c$ and suppose $f\le 1$ and
$$\nu = \cases{\bigl\lfloor \sqrt{2/f} \bigr\rfloor
      ,&if $2\le \bigl\lfloor \sqrt{2/f} \bigr\rfloor$;\cr
  2,&if $\bigl\lfloor \sqrt{2/f}\bigr\rfloor < 2$.\cr}$$
Then a single ring with $\nu$~nodes and capacity~$c$
can support uniform traffic of $d$~units.
\endlem

\prf{Proof}
Except for the fact that traffic can be split,
this is essentially a ring loading problem,
so the reader may expect us to apply known results for ring loading.
However, a na\"{\i}ve method suffices for our present purposes:
For each pair of nodes,
route $\lfloor d/2 \rfloor$ units of traffic
on the major (longer) arc
and $\lceil d/2\rceil$ units of traffic on the minor arc.
When both arcs are of equal length,
either one may be designated the major arc.
There are four cases:
$d$ may be even ($d=2\delta$) or odd ($d=2\delta+1$),
and independently, $\nu$ may be even ($\nu=2\mu$) or odd ($\nu=2\mu+1$).
In three of these cases it is easy to calculate
the amount of traffic on the most heavily loaded link.
The exceptional case is $\nu=2\mu$ and $d=2\delta+1$,
where the maximum load depends on
exactly how the major/minor arcs are chosen
for diametrically opposite pairs of vertices.
In this case we obtain an upper bound by observing
that the worst case occurs when
there is some edge of~$C_n$ that is systematically chosen
to be the ``major arc'' in all these tie-breaking situations.

The table below summarizes the results of these calculations.
$$\vbox{\offinterlineskip
\hrule
\halign{&\vrule#&\strut\quad#\hfil\quad\cr
height2pt&\omit&&\omit&&\omit&\cr
&\omit&&\hfil$d=2\delta$&&\hfil$d=2\delta+1$&\cr
height2pt&\omit&&\omit&&\omit&\cr
\noalign{\hrule}
height2pt&\omit&&\omit&&\omit&\cr
&$\nu=2\mu$&&$\mu(2\mu-1)\delta$
  &&\hfil${1\over2}\mu\bigl[(\mu+1)(\delta+1) + (3\mu-3)\delta\bigr]$&\cr
height2pt&\omit&&\omit&&\omit&\cr
\noalign{\hrule}
height2pt&\omit&&\omit&&\omit&\cr
&$\nu=2\mu+1$&&$\mu(2\mu+1)\delta$
  &&\hfil${1\over2}\mu\bigl[(\mu+1)(\delta+1) + (3\mu+1)\delta\bigr]$&\cr
height2pt&\omit&&\omit&&\omit&\cr}
\hrule}$$
(In the right-hand column,
the two summands correspond to the minor arcs and major arcs
respectively that contain the most heavily loaded link.)
We need to prove that in every case,
$c$ is large enough to accommodate the most heavily loaded link.
First let us suppose that $\nu=2$.
Then the first row of the table applies, with $\mu=1$;
the two entries reduce to $\delta$ and $\delta+1$ respectively,
and in either case this is at most $\lceil d/2\rceil$.
Since $f\le 1$ and $c$ is an integer, we have $c\ge \lceil d/2\rceil$,
which shows that there is enough capacity if $\nu=2$.

If $\nu\ne 2$, then we must have $\nu \le \sqrt{2/f}$,
which implies that $c \ge \nu^2d/4$.
To check that there is sufficient capacity,
we therefore just need to check,
for each of the four entries in the above table,
that if we subtract it from $\nu^2d/4$, the result is always nonnegative.
If $\nu=2\mu$ and $d=2\delta+1$ then
$${\nu^2d\over 4} - 
  {\mu\bigl[(\mu+1)(\delta+1) + (3\mu-3)\delta\bigr] \over 2}
 = {\mu(\mu+2\delta -1)\over 2},$$
which is always nonnegative (since, for example, $\mu\ge1$).
In the other three cases, it is easily checked that the analogous calculation
yields a polynomial in $\mu$ and~$\delta$ with nonnegative coefficients,
which is therefore always nonnegative.\qquad\endproof

Finally, we need to show that Algorithm~A
is indeed an approximation algorithm.

\thm{Theorem 2}
Algorithm~A is a $12\sqrt2$-approximation algorithm.
\endthm

\prf{Proof}
If $f\ge1$, then Algorithm~A uses $n(n-1)\lceil f\rceil$ ADM's,
and the LP bound is $n(n-1)f$.  Dividing the former by the latter
yields $\lceil f\rceil/f$, which is at most~$2$ because $f\ge 1$.
So this case is all right.

Now suppose $f<1$.  The number of ADM's used is at most
$$m = {\lceil n/\mu\rceil \choose 2}M,$$
where $\mu = \lfloor M/2\rfloor$.
We need to show that $m$ is
at most $12\sqrt2$ times the minimum number of ADM's.
We begin by noting that the case $M=n$ is easy:
Algorithm~A generates a single ring with $n$~ADM's,
which Lemma~1 tells us can handle all the demands,
and since there is some traffic terminating at every vertex,
this solution is optimal.
Therefore we assume for the rest of the proof that $M<n$,
which implies in particular that $\sqrt{2/f} < M+1$
(by definition of~$M$).

Let $q$ denote $m$ divided by the lower bound of Theorem~1.
Consider first the case $\mu=1$.
If $\mu=1$ then $M\le 3$ and $\sqrt{2/f} < 4$ so
$$q = {m \over {1\over4}(n^2-1)\sqrt{f}} = {2n(n-1)M \over(n+1)(n-1)\sqrt{f}}
  < {12n\sqrt{2} \over n+1} < 12\sqrt{2},$$
as required.

So suppose that $\mu\ge2$.
Now $\lceil n/\mu\rceil \le (n+\mu-1)/\mu$ and $M\le 2\mu+1$ so
$$m \le {1\over2} \biggl({n+\mu-1\over \mu}\biggr)
  \biggl({n-1 \over \mu}\biggr) (2\mu+1).$$
Dividing by ${1\over 4}(n^2-1)\sqrt{f}$
and using the inequality $\sqrt{2/f} < M+1 \le 2\mu+2$ yields
$$q \le 2\biggl({n+\mu-1\over n+1}\biggr) {2\mu+1\over \mu^2\sqrt{f}}
  \le 2\sqrt{2}\biggl({n+\mu-1\over n+1}\biggr) {(2\mu+1)(\mu+1)\over \mu^2}.$$
Now $\mu-1 < M/2 < n/2$ so
$${n+\mu-1\over n+1} < {3n/2\over n+1} < {3\over 2}.$$
Also $(2\mu+1)(\mu+1)/\mu^2$ is a decreasing function of~$\mu$
so it attains its maximum when $\mu$ is as small as possible,
i.e., when $\mu=2$.
Therefore
$$q < 2\sqrt{2} \biggl({3\over 2}\biggr){(2\cdot 2 + 1)(2 + 1)\over 2^2}
   = {45\sqrt{2}\over 4} < 12\sqrt2,$$
completing the proof.
\qquad\endproof

\cor{Corollary 2}
There is an $\alpha$-approximation algorithm for ring grooming
with $K$-quasi-uniform traffic, where $\alpha = \max(2K, 12\sqrt{2K})$.
\endcor

\prf{Proof}
Given an instance~$I$ of ring grooming with $K$-quasi-uniform traffic,
we let $d$ be the size of the largest $d_{jk}$, and we
we create an instance~$I'$ of ring grooming with uniform traffic
by changing every $d_{jk}$ (with $j\ne k$) to~$d$.
Then we apply Algorithm~A to~$I'$ to obtain a solution
that a~fortiori yields a solution to~$I$.
Now we trace through the proof of Theorem~2.
If $d\ge 2c$ then at most $n(n-1)\lceil d/2c\rceil$ ADM's are used,
but between any two vertices we have at least $d/K$ units of traffic,
yielding a lower bound of $n(n-1)d/(2cK)$.
Dividing the former by the latter
gives a ratio of at most~$2K$, since $d\ge 2c$.

Otherwise, $d<2c$, and the argument in the proof of Theorem~2
for the case $M=n$ still applies to show that
the solution is optimal in this case.
If $d<2c$ and $M<n$, then
let $q'$ be $m$ divided by the lower bound of Corollary~1.
Then $q'=q\sqrt{K}$, where $q$ is defined
in the proof of Theorem~2.
That proof shows that
$q<12\sqrt{2}$, whence $q'<12\sqrt{2K}$.\qquad\endproof

Although we are not proposing that Algorithm~A
be used in a practical implementation,
it is of some interest to consider whether
the constant~$12\sqrt2$ can be improved.
One possibility is to use a more sophisticated method
for constructing $(n,M,2)$-covering designs (see for example~[10]).
For a non-trivial example of what a good design can accomplish,
suppose that $n=15$ and $c=d=1$.
The add/drop lower bound is~$105$.
Our approximation algorithm tells us to take $M=2$,
but suppose we take $M=3$ instead
and look for a $(15,3,2)$-covering design.
One example of this is the solution to
the famous {\it Kirkman schoolgirl problem\/}
(see [3] or any standard reference on combinatorial designs).
This yields an optimal solution with 35~rings
and three ADM's per ring.

In general, however, we cannot expect to find such excellent designs fast.
Moreover, a trivial lower bound on the size of an $(n,M,2)$-covering design
is ${n \choose 2}/{M \choose 2}$,
and if we carry this through the proof of Theorem~2,
we find that the maximum possible improvement
in the guaranteed performance ratio
is a factor of two or three.

Another possible improvement is
to use a better ring loading method in Lemma~1.
By considering a maximum cut, one can show that
the largest improvement factor we can hope for here is~$\sqrt2$.
Achieving this is not easy, though;
for example, it is certainly not possible to na\"{\i}vely replace
$\sqrt2$ with~$2$ in Lemma~1:
If $c=2$ and $d=1$ then $M$ can be at most three,
whereas $\bigl\lfloor 2/\sqrt{f}\bigr\rfloor = 4$.

Finally, perhaps the bound of Theorem~1 can be improved.
For the special case of $d=1$ we have the following argument
that gives a slightly better lower bound.
The idea is that at any particular ADM,
the $t$ units of traffic exiting from a particular side of the ADM
all have different destinations, say $v_1, v_2, \ldots, v_t$.
Because the network is a ring,
this means that the traffic destined for~$v_j$
must pass through the ADM's at $v_i$ for all $i<j$,
thus ``wasting'' some of the capacity of these ADM's.
If we work through the details of this
for the example given at the end of section~5,
we obtain a lower bound of about ${1\over2}n^{3/2}$,
a factor-of-two improvement.
Unfortunately, we do not see how to generalize this argument to arbitrary~$d$.

\heading{7}{Open problems and related questions}

Of course, the main open problem is to find
a constant-factor approximation algorithm
for the general ring grooming problem,
or to prove that one of the algorithms in the literature
is in fact a constant-factor approximation algorithm.
As we mentioned above, it is also open whether
ring grooming is fixed-parameter tractable when $n$ is the parameter.
It would be of practical interest if the answer turns out to be yes,
since SONET standards do not permit more than sixteen ADM's on a ring.

If a good approximation algorithm cannot be found,
then an inapproximability result would be desirable.
We suspect that ring grooming is MAXSNP-hard but we cannot prove it.
Ring grooming is reminiscent of several other well-known problems,
e.g., integer multicommodity flow, clique cover, and facility location,
but we have been unable to construct an explicit reduction.
Note that comparing Theorem~1 and the LP bound
reveals a large integrality gap.
We expect that our ILP formulation can be improved;
for example, Michel Goemans (personal communication) has pointed out
that the following constraints can be added.
$$\leqalignno{
 t^0_{ijk} + t^1_{ijk}
  &\le \hbox to 1in{$d_{jk}x_{ij}$,\hfil} \forall i, j, k\; (j<k)&(5a)\cr
 t^0_{ijk} + t^1_{ijk}
  &\le \hbox to 1in{$d_{jk}x_{ik}$,\hfil} \forall i, j, k\; (j<k)&(5b)\cr
}$$
Summing these inequalities over~$i$ and using constraint~(1),
we deduce that $\sum_i x_{ij} \ge 1$ whenever
there is any traffic terminating at vertex~$j$,
even if the $d_{jk}$ are all small and
the $x_{ij}$ are allowed to be arbitrary real numbers.
While this is an improvement, further constraints would be desirable.
Unfortunately, they seem to be hard to find.

Our formulation of the ring grooming problem is not the only possible one;
other authors have considered slightly different versions.
We describe here some of the more important variations.

\meti{(1)} {\it Unidirectional rings}.  Unlike ring loading,
ring grooming is nontrivial even for unidirectional rings,
and as mentioned in the introduction,
it has been studied by some of the authors listed in the references.
Note that some service providers tend to use unidirectional rings
in situations where there exists a {\it hub node\/}
at which all traffic terminates,
and standard bin packing algorithms work well here.
\meti{(2)} {\it No timeslot interchange}.
In our formulation, we regard a routing as feasible
provided that the total number of units of traffic on each edge of~$C_n$
does not exceed~$c$.
In actual SONET rings,
the $c$ units of capacity are $c$ separate {\it timeslots}.
We have assumed that all our ADM's have {\it timeslot interchange\/} capability,
meaning that whenever a unit of traffic passes through an ADM,
it can be switched to a different timeslot if necessary.
Some less expensive ADM's do not have this capability,
in which case all traffic must choose a timeslot
and remain on the same timeslot from source to destination.
This creates some additional constraints that are akin to graph coloring.
\meti{(3)} {\it Cross-connection between rings}.
In our formulation, we require that
traffic stay on the same ring from source to destination.
In some networks, there are {\it digital cross-connects\/}
installed in the central offices which allow traffic to switch from
one ring to another at a vertex provided that
both rings have an ADM at that location.
If we assume that cross-connection is available,
then ring grooming becomes more similar to multicommodity flow,
and perhaps some of the techniques carry over.
\meti{(4)} {\it No traffic splitting}.
As we have already remarked,
our formulation allows traffic to be split between major and minor arcs,
as long as all variables remain integers.
One might wonder whether this assumption is realistic.
Certainly if $d_{jk} > c$ then $d_{jk}$ must at least
be split across different rings,
and if splitting across rings is permissible then
splitting between major and minor arcs is likely to be permissible too.
However, what if $d_{jk} \le c$?
Typically, in real networks,
traffic demands come in multiples of a few standard units,
e.g., the traffic between a particular pair of vertices
might be a combination of several units of size twelve
and several units of size three.
It would not be permissible to split a size-three or size-twelve unit,
either across rings or between major and minor arcs,
but for example the various units of size three could be routed independently.
It would be interesting to study the effect of multiple unit sizes,
not only on the ring grooming problem, but also on the ring loading problem.
\meti{(5)} {\it Unidirectional traffic and asymmetric routing}.
We have assumed not only that our rings are bidirectional,
but that our traffic is bidirectional as well,
and moreover that the traffic from $j$ to~$k$
must be routed the same way as the traffic from $k$ to~$j$,
only in reverse.
These are reasonable assumptions for current networks,
but some experts predict future increases
in unidirectional traffic ($d_{jk} \ne d_{kj}$).
Also, even if $d_{jk} = d_{kj}$,
it is in principle possible to route the two directions independently.
\meti{(6)} {\it Dynamic traffic}.
In real networks the traffic matrix is usually not fixed but grows with time.
Many network operators have a {\it cap-and-grow\/} policy,
meaning that existing traffic is not allowed to be rerouted when new traffic arrives.
There is therefore a practical need for an {\it online algorithm\/}
for ring grooming.  Very little work has been done so far on this
important but difficult problem.

\heading{8}{Acknowledgments}
The first author would like to thank Xingxing Yu and Michel Goemans
for useful conversations and Oren Patashnik for pointing him
in the right direction regarding covering designs.
Most of the work in this paper was done
while the authors were at the Tellabs Research Center in Cambridge.

\noindent
\Refs
 
\ref 1\\
{\smc N. Brauner, Y. Crama, G. Finke, P. Lemaire, and C. Wynants},
{\it Approximation algorithms for SONET networks},
Technical Report GEMME 2002, University of Li\`ege, 2002.
\endref

\ref 2\\
{\smc S. Cosares and I. Saniee},
{\it An optimization problem related to balancing loads on SONET rings},
Telecom.\ Systems, 3 (1994), pp.~165--181.
\endref

\ref 3\\
{\smc H. D\"orrie},
{\it 100 Great Problems of Elementary Mathematics:
Their History and Solution},
Dover, New York, 1965.
\endref

\ref 4\\
{\smc R. G. Downey and M. R. Fellows},
{\it Parameterized Complexity},
Springer, New York, 1999.
\endref

\ref 5\\
{\smc A. Frank, T. Nishizeki, N. Saito, H. Suzuki, and E. Tardos},
{\it Algorithms for routing around a rectangle},
Discrete Appl.\ Math., 40 (1992), pp.~363--378.
\endref

\ref 6\\
{\smc M. R. Garey and D. S. Johnson},
{\it Computers and Intractability:
A Guide to the Theory of NP-Completeness},
W. H. Freeman, San Francisco, 1979.
\endref

\ref 7\\
{\smc O. Gerstel, P. Lin, and G. Sasaki},
{\it Combined WDM and SONET network design},
Proc.\ IEEE INFOCOM '99, pp.~734--743.
\endref

\ref 8\\
{\smc O. Gerstel, P. Lin, and G. Sasaki},
{\it Wavelength assignment in a WDM ring
to minimize cost of embedded SONET rings},
Proc.\ IEEE INFOCOM '98, pp.~94--101.
\endref

\ref 9\\
{\smc O. Goldschmidt, D. S. Hochbaum, A. Levin, and E. V. Olinick},
{\it The SONET edge-partition problem},
Networks, 41 (2002), pp.~13--23.
\endref

\ref 10\\
{\smc D. M. Gordon, O. Patashnik, and G. Kuperberg},
{\it New constructions for covering designs},
J. Combin.\ Designs, 3 (1995), pp.~269--284;
http://sdcc12.ucsd.edu/\~{ }xm3dg/errata.html
\endref

\ref 11\\
{\smc P. Harshavardhana, P. K. Johri, and R. Nagarajan},
{\it A note on weight-based load balancing on SONET rings},
Telecom.\ Systems, 6 (1996), pp.~237--239.
\endref

\ref 12\\
{\smc J-Q Hu},
{\it Traffic grooming in wavelength-division-multiplexing ring networks:
A linear programming solution},
J. Optical Networking, 1 (2002), pp.~397--408.
\endref

\ref 13\\
{\smc S. Khanna},
{\it A polynomial time approximation scheme for the SONET ring
loading problem}, Bell Labs Tech.\ J.,
%
%
2 (Spring 1997), pp.~36--41.
\endref

\ref 14\\
{\smc Y. Lee, H. Sherali, J. Han, and S. Kim},
{\it A branch-and-cut algorithm for solving
an intra-ring synchronous optical network design problem},
Networks, 35 (2000), pp.~223--232.
\endref

\ref 15\\
{\smc X.-Y. Li, L.-W. Liu, P.-J. Wan, and O. Frieder},
{\it Practical traffic grooming scheme for single-hub SONET/WDM rings},
Proc.\ IEEE LCN '00, pp.~556--564.
\endref

\ref 16\\
{\smc E. H. Modiano and A. L. Chiu},
{\it Traffic grooming algorithms for reducing electronic multiplexing costs
in WDM ring networks},
J. Lightwave Tech., 18 (2000), pp.~2--12.
\endref

\ref 17\\
{\smc Y.-S. Myung, H.-G. Kim, and D.-W. Tcha},
{\it Optimal load balancing on SONET bidirectional rings},
Operations Research, 45 (1997), pp.~148--152.
\endref
 
\ref 18\\
{\smc A. Schrijver, P. Seymour, and P. Winkler},
{\it The ring loading problem},
SIAM J. Discrete Math., 11 (1998), pp.~1--14.
\endref

\ref 19\\
{\smc J. M. Simmons, E. L. Goldstein, and A. A. M. Saleh},
{\it Quantifying the benefit of wavelength add-drop in WDM rings
with distance-independent and dependent traffic},
J. Lightwave Tech., 17 (1999), pp.~48--57.
\endref

\ref 20\\
{\smc A. Sutter, F. Vanderbeck, and L. Wolsey},
{\it  Optimal placement of add/drop multiplexers:
Heuristic and exact algorithms},
Operations Research, 46 (1998), pp.~719--728.
\endref

\ref 21\\
{\smc P.-J. Wan, L.-W. Liu, and O. Frieder},
{\it Grooming of arbitrary traffic in SONET/WDM BLSRs},
IEEE J. Selected Areas Comm., 18 (2000), pp.~1995--2003.
\endref

\ref 22\\
{\smc J. Wang, W. Cho, V. Vemuri, and B. Mukherjee},
{\it Improved approaches for cost-effective traffic grooming
in WDM ring networks:
non-uniform traffic and bidirectional ring},
Proc.\ IEEE ICC '00, pp.~1295--1299.
\endref

\ref 23\\
{\smc X. Zhang and C. Qiao},
{\it An effective and comprehensive approach
for traffic grooming and wavelength assignment in SONET/WDM rings},
IEEE/ACM Trans.\ Networking, 8 (2000), pp.~608--617.
\endref

\bye